\documentclass[a4paper,11pt] {amsart}
\usepackage{amssymb,latexsym,amsmath}
\usepackage{amsfonts,amsbsy,bm}
\usepackage{color}
\usepackage{graphicx}
\usepackage{listings}
\usepackage{float}
\usepackage{multirow}
\usepackage{graphicx}
\usepackage{booktabs}
\usepackage{varwidth}
\providecommand{\tabularnewline}{\\}
\newcommand{\bprop} {\begin{proposition}}
\newcommand{\eprop} {\end{proposition}}
\newcommand{\btheo} {\begin{theorem}}
\newcommand{\etheo} {\end{theorem}}
\newcommand{\blem} {\begin{lemma}}
\newcommand{\elem} {\end{lemma}}
\newcommand{\bcor} {\begin{corollary}}
\newcommand{\ecor} {\end{corollary}}

\newcommand{\Be}{\begin{equation}}
\newcommand{\Ee}{\end{equation}}
\newcommand{\Bea}{\begin{eqnarray}}
\newcommand{\Eea}{\end{eqnarray}}
\newcommand{\Bes}{\begin{equation*}}
\newcommand{\Ees}{\end{equation*}}
\newcommand{\Beas}{\begin{eqnarray*}}
\newcommand{\Eeas}{\end{eqnarray*}}
\newcommand{\Ba}{\begin{array}}
\newcommand{\Ea}{\end{array}}

\definecolor{Brown}{cmyk}{0,0.81,1,0.60}
\definecolor{OliveGreen}{cmyk}{0.64,0,0.95,0.40}
\definecolor{CadetBlue}{cmyk}{0.62,0.57,0.23,0}
\definecolor{lightlightgray}{gray}{0.91}

\begin{document}
\lstset{
	language=Python,                             % Code langugage
	basicstyle=\ttfamily,                   % Code font, Examples: \footnotesize, \ttfamily
	keywordstyle=\color{black},        % Keywords font ('*' = uppercase)
	commentstyle=\color{gray},              % Comments font
	numbers=left,                           % Line nums position
	numberstyle=\tiny,                      % Line-numbers fonts
	stepnumber=1,                           % Step between two line-numbers
	numbersep=5pt,                          % How far are line-numbers from code
	backgroundcolor=\color{lightlightgray}, % Choose background color
	frame=none,                             % A frame around the code
	tabsize=2,                              % Default tab size
	captionpos=b,                           % Caption-position = bottom
	breaklines=true,                        % Automatic line breaking?
	breakatwhitespace=false,                % Automatic breaks only at whitespace?
	showspaces=false,                       % Dont make spaces visible
	showtabs=false,                         % Dont make tabls visible
	morekeywords={__global__, __device__},  % CUDA specific keywords
}
\theoremstyle{plain}% Theorem-like structures
\newtheorem{theorem}{Theorem}[section]
\newtheorem{corollary}[theorem]{Corollary}
\newtheorem{lemma}[theorem]{Lemma}
\newtheorem{proposition}[theorem]{Proposition}

\theoremstyle{definition}
\newtheorem{definition}[theorem]{Definition}
\newtheorem{example}[theorem]{Example}

\theoremstyle{remark}
\newtheorem{remark}[theorem]{Remark}
%%%%%%%%%%%%%%%%%%%%%%%%%%%%%%%%%%%%%%%%%%%%%%%%%%%%%%%%%%%%%%%%%%%%
\title[Some inequalities for the beta function and its ratios]{Some inequalities for the beta function and its ratios }
%--------author
%\author{Noname}
\author[J.M. Tanoh Dje]{Jean$-$Marcel T. Dje}
\address{Unit\'e de Recherche et d'Expertise Num\'erique, Universit\'e Virtuelle de C\^ote d'Ivoire, Cocody II-Plateaux - 28 BP 536 ABIDJAN 28}
\email{{\tt tanoh.dje@uvci.edu.ci}}

%\author[A. Asare-Tuah]{Anton Asare-Tuah}
%\address{Department of Mathematics, University of Ghana,\\ P. O. Box LG 62 Legon, Accra, Ghana}
%\email{aasare-tuah@ug.edu.gh}

%\author[E. Djabang]{Emmanuel Djabang}
%\address{Department of Mathematics, University of Ghana,\\ P. O. Box LG 62 Legon, Accra, Ghana}
%\email{EDjabang@ug.edu.gh}

\author[E. A. K. Schwinger]{Eyram A. K. Schwinger}
\address{Department of Mathematics, University of Ghana,\\ P. O. Box LG 62 Legon, Accra, Ghana}
\email{eakschwinger@ug.edu.gh}

 \author[B. F. Sehba]{Benoit F. Sehba}
\address{Department of Mathematics, University of Ghana,\\ P. O. Box LG 62 Legon, Accra, Ghana}
\email{bfsehba@ug.edu.gh}

%\author[R. A. Twum]{Ralph A. Twum}
%\address{Department of Mathematics, University of Ghana,\\ P. O. Box LG 62 Legon, Accra, Ghana}
%\email{ratwum@ug.edu.gh}


\maketitle

%-------abstract
\begin{abstract}
In this paper, we prove some inequalities for the differences and ratios of the beta function.
\end{abstract}
%  \textit{This paper is dedicated to Dr. Douglas Adu-Gyamfi in honor of his 80th birthday.}
%%% ----------------------------------------------------------------------
%\maketitle
%%% ----------------------------------------------------------------------
%\vspace*{\fill} % Pushes the text to the vertical center
%\vspace*{\fill}
\section{Introduction}

Recall that the beta function is defined by 
\Be
B(x,y):=\int_0^1t^{x-1}(1-t)^{y-1}dt,\quad x>0, y>0.
\Ee
A related function is the gamma function defined by
$$\Gamma(x):=\int_0^\infty t^{x-1}e^{-t}dt,\quad x>0.$$
Indeed, it is not difficult to check that
$$B(x,y)=\frac{\Gamma(x)\Gamma(y)}{\Gamma(x+y)}.$$
For more information on the above functions and related functions, see \cite{Abramowitz}.
\medskip

An interesting question on the beta function is the problem of finding 'good or best' functions $f$ and $g$ such that
\Be\label{eq:problem1}
f(x,y)\leq B(x,y)\leq g(x,y)
\Ee
and recently,
\Be\label{eq:problem2}
f(a,b,y_1,y_2)\leq \frac{B(b,y_2)}{B(a,y_1)}\leq g(a,b,y_1,y_2).
\Ee
For the first problem (\ref{eq:problem1}), one of the first noticeable results is the inequality described by S. Dragomir et al. \cite[Theorem 3]{Dragomir}:
\Be\label{eq:dragomir}
\frac{1}{xy}-\frac{1}{4}\leq B(x,y)\leq \frac{1}{xy}, x>1, y>1.
\Ee
The above inequality has since been improved and extended by various authors \cite{Alzer1,Alzer,Alzer2,Cerone,From,Grenie,Ivady}. In this direction, let us mention the following estimates obtained by P. Ivady in \cite{Ivady}.
\Be\label{eq:ivady}
\frac{x+y-xy}{xy}\leq B(x,y)\leq\frac{x+y}{xy(1+xy)},\quad\text{for}\quad 0<x,y\leq 1.
\Ee

Several one-sided (lower or upper) inequalities were obtained in \cite{Cerone,From,Grenie}. The following will be relevant in our discussion.
% $0<x,y<\infty$ with $x>2$ and $xy<1$, compare
%$\frac{x^{x-1}y^{y-1}}{(x+y)^{x+y-1}}$
%$\frac{(x+y)(1-\sqrt{xy})^2}{xy(x+1)}(y+1)$ and

%$\frac{x^xy^y}{(x+y)^{x+y}}\sqrt{\frac{2\pi}{x}}$. 
\btheo\label{thm:17}\cite[Page 2]{Grenie}
For $x>0$,$y>0$,
$$B(x,y)\geq \frac{x^{x-1}y^{y-1}}{(x+y)^{x+y-1}}.$$
\etheo
A second result that we will be comparing with our findings is given as follows.
\btheo\label{thm:31}\cite[Theorem 3.1]{From}
Let $x>0$,$y>0$ satisfy $xy<1$. Then
 $$B(x,y)\geq \frac{(x+y)(1-\sqrt{xy})^2}{xy(x+1)(y+1)}.$$
\etheo
The following left-side inequality is also from \cite{From}.
\btheo\label{thm:32}\cite[Theorem 3.2]{From}
Suppose that $x>0$, $y>0$ with $\frac{x}{x+1}+\frac{y}{y+1}<1$. Then
 
$$B(x,y)\geq \frac{x+y}{x}\cdot\frac{z^2}{(z^2+1)y}$$ where
$$z=\frac{1-\left(\frac{x}{x+1}+\frac{y}{y+1}\right)}{\sqrt{\frac{x}{(x+1)^2(x+2)}+\frac{y}{(y+1)^2(y+2)}}}.$$
\etheo
We will also compare our finding with the following non-trivial result.
\btheo\label{thm:34}\cite[Theorem 3.4]{From}
Suppose that $x>0$, $y>0$ with $xy<1$. Then
 
$$B(x,y)\geq \frac{x+y}{xy}[1-u(x,y)]$$ where
$$u(x,y)=e^{-s}\left[1+s\left(\frac{x}{x+1}\right)+(e^s-1-s)\left(\frac{x}{x+2}\right)\right]\cdot\left[1+s\left(\frac{y}{y+1}\right)+(e^s-1-s)\left(\frac{y}{y+2}\right)\right]$$
and $$s=\frac{1}{2}\ln\left(\frac{1}{xy}\right).$$
\etheo
The following is the first part of \cite[Theorem 3.11]{From}.
\btheo\label{thm:311}
For $x\geq 3$, and $y>0$,
 
$$B(x,y)\geq \frac{1}{y}\left[\left(\frac{1}{y+1}\right)^{x-1}+\frac{\frac{1}{2}(x-1)y}{(y+1)^{x-1}(y+2)}\right].$$ 
\etheo
The last lower estimate we recall is again from \cite{From}.
\btheo\label{thm:314}\cite[Theorem 3.14]{From}
Let $x>0$, and $y>0$. Then
 
$$B(x,y)\geq \left(\frac{x}{x+y}\right)^xe^{-x}\sqrt{\frac{2\pi}{x}}.$$ 
\etheo
\medskip

The problem (\ref{eq:problem2}) was recently considered in \cite[Theorem 3.12, 3.16]{From} where several estimates were proved. Let us mention two of these inequalities:
\btheo\label{thm:From}
Let $0<a<b$ and $y>0$. Then
\Be\label{eq:Fromlowerupper}
\frac{b^{b-1}(a+y)^{a+y-1}}{a^{a-1}(b+y)^{b+y-1}}\leq \frac{B(b,y)}{B(a,y)}\leq \frac{b^{b}(a+y)^{a+y}}{a^{a}(b+y)^{b+y}}.
\Ee
\etheo
Let us mention the following result recently proved in \cite[Theorem 4.12]{UG} as follows.
\btheo\label{thm:twoside1}
For $a,b,y>0$ with $a<b$, the following holds.
\Be\label{eq:twoside1}
\frac{b+y}{a+y}\left(\frac{a}{b}\right)^{y+1}<\frac{B(b,y)}{B(a,y)}<\frac{b+y}{a+y}\left(\frac{a}{b}\right)\left(\frac{a+y+1}{b+y+1}\right)^y.
\Ee
\etheo
We recall that the lower bound in (\ref{eq:Fromlowerupper}) is better than the one in (\ref{eq:twoside1}) but the inequality between the two upper bounds depends on values of $a,b$ and $y$.
\medskip

The following estimate was also obtained in \cite[Theorem 3.17]{From}.
\btheo\label{thm:From317}
Assume that $0<a<b<\infty$, $y>0$. Let 
$$c=1-\frac{y}{b-a}\ln\left(\frac{b+y}{a+y}\right),\quad \alpha=\ln c-1,\quad \beta=\frac 1c.$$
Then
\Be\label{eq:From317}
\frac{B(b,y)}{B(a,y)}\leq e^{(\alpha+\beta)(b-a)}\left(\frac{a+y}{b+y}\right)^{y\beta}.
\Ee
\etheo
%\Be\label{eq:Fromlower}
%\frac{B(b,y)}{B(a,y)}\geq \frac{b^{b-1}(a+y)^{a+y-1}}{a^{a-1}(b+y)^{b+y-1}},\quad 0<a<b, y>0.
%\Ee
%\Be\label{eq:Fromupper}
%\frac{B(b,y)}{B(a,y)}\leq \frac{b^{b}(a+y)^{a+y}}{a^{a}(b+y)^{b+y}},\quad 0<a<b, y>0.
%\Ee
In this paper, our objective is to prove new inequalities of the form (\ref{eq:problem2}) that improve (\ref{eq:Fromlowerupper}), (\ref{eq:twoside1}), and (\ref{eq:From317}), and estimates of the differences of the beta function. By the latter, we mean inequalities of the form
\Be\label{eq:problem3}
f(a,b,y)\leq B(b,y)-B(a,y)\leq g(a,b,y),\quad 0<a<b, y>0.
\Ee
Note that (\ref{eq:problem1}) can be obtained from (\ref{eq:problem3}) for some parameters $a,b,y$. 
\medskip

We find ranges of $a,b,y$ where our results improve (\ref{eq:Fromlowerupper}) and (\ref{eq:twoside1}). We prove this manually, but we also provide a table that compares our bounds with those in Theorem \ref{thm:From} and Theorem \ref{thm:twoside1}. We also obtain some inequalities for the differences of the beta function. This result then leads to a lower bound for the beta function that we will compare with Theorem \ref{thm:17}, \ref{thm:31}, \ref{thm:32}, \ref{thm:34}, \ref{thm:311} and \ref{thm:314}. For the analog of (\ref{eq:From317}), using Python, we show in a figure the area of the parameters where our result is an improvement. The associated Python code is provided in Section 7. Our results are presented and proved in the next section.
%\section{Preliminaries}
\section{Main results}
A key lemma in this work is the following.
\blem[Generating lemma]\label{lem:keylem}
Let $k$ be a positive integer and let $h$ and $l$ be real numbers. Then for $x,y>0$, the following identity holds.
\Be\label{eq:integparts1}\int_0^1 t^{x+k-l-1}(1-t)^{y-h} dt+(x+k-l)\int_0^1 t^{x+k-l-1}(1-t)^{y-h}\ln t dt=(y-h)\int_0^1 t^{x+k-l}(1-t)^{y-h-1}\ln t dt.
\Ee
\elem

\begin{proof}
Put $$I=\int_0^1 t^{x-l}(1-t)^{y-h}(t^k\ln t)' dt.$$  
We compute $I$ in two different ways. First, we have
\Be\label{eq:integparts11}
I =\int_0^1 t^{x-l}(1-t)^{y-h}(kt^{k-1}\ln t+t^{k-1}) dt.
\Ee
Second, using integration by parts, we obtain
\Be\label{eq:integparts2}
I = -(x-l)\int_0^1 t^{x+k-l-1}(1-t)^{y-h}\ln t dt+(y-h)\int_0^1 t^{x+k-l}(1-t)^{y-h-1}\ln t dt.
\Ee
It follows from (\ref{eq:integparts11}) and (\ref{eq:integparts2}) that 
\Beas
\int_0^1 t^{x-l}(1-t)^{y-h}(kt^{k-1}\ln t+t^{k-1}) dt &=& -(x-l)\int_0^1 t^{x+k-l-1}(1-t)^{y-h}\ln t dt\\ &+& (y-h)\int_0^1 t^{x+k-l}(1-t)^{y-h-1}\ln t dt. 
\Eeas
Rearranging the terms, we then obtain
$$\int_0^1 t^{x+k-l-1}(1-t)^{y-h} dt+(x+k-l)\int_0^1 t^{x+k-l-1}(1-t)^{y-h}\ln t dt=(y-h)\int_0^1 t^{x+k-l}(1-t)^{y-h-1}\ln t dt.$$
\end{proof}

We have the following first result.
\btheo\label{thm:main1}
Let $0<a<b$. The following hold.
\begin{itemize}
    \item[(1)] If $0<y\leq 1$, then
    \Be\label{eq:main11}
    \left(\frac{b+y}{a+y}\right)^{1-y}\left(\frac{a}{b}\right)\leq\frac{B(b,y)}{B(a,y)}\leq \left(\frac{a}{b}\right)^{y}.
    \Ee
    \item[(2)] If $y>1$, then 
    \Be\label{eq:main12}
    \left(\frac{a}{b}\right)^{y}\leq\frac{B(b,y)}{B(a,y)}\leq \left(\frac{b+y}{a+y}\right)^{1-y}\left(\frac{a}{b}\right).
    \Ee
\end{itemize}
\etheo
\begin{proof}
    Taking $k=l$ and $h=1$ in (\ref{eq:integparts1}), we obtain
    \Be\label{eq:integparts3}\int_0^1 t^{x-1}(1-t)^{y-1} dt+x\int_0^1 t^{x-1}(1-t)^{y-1}\ln t dt=(y-1)\int_0^1 t^{x}(1-t)^{y-2}\ln t dt.
\Ee
We assume that $0<y\leq 1$. Then we read the left hand side of (\ref{eq:integparts3}) as 
$$\int_0^1 t^{x-1}(1-t)^{y-1} dt+x\int_0^1 t^{x-1}(1-t)^{y-1}\ln t dt=B(x,y)+x\partial_xB(x,y).$$
Let $$I=(y-1)\int_0^1 t^{x}(1-t)^{y-2}\ln t dt.$$
We us recall that for $\frac{(t-1)}{t}\leq \ln t\leq t-1$ when $0<t\leq 1$. It follows that
$$-(y-1)\int_0^1 t^{x}(1-t)^{y-1} dt\leq I\leq -(y-1)\int_0^1 t^{x-1}(1-t)^{y-1} dt.$$
We deduce that
\Bea
-(y-1)B(x+1,y)\leq B(x,y)+x\partial_xB(x,y)\leq -(y-1)B(x,y)
\Eea
or equivalently,
\Be\label{eq:diffquotient1}
-\frac{y-1}{x+y}-\frac{1}{x}\leq \frac{\frac{\partial B(x,y)}{\partial x}}{B(x,y)}\leq -\frac{y}{x}.
\Ee
Integrating the above equation from $x=a$ to $x=b$, we obtain
$$-(y-1)\ln\left(\frac{b+y}{a+y}\right)-\ln\left(\frac{b}{a}\right)\leq \ln\left(\frac{B(b,y)}{B(a,y)}\right)\leq -y\ln\left(\frac{b}{a}\right).$$
That is 
$$ \left(\frac{b+y}{a+y}\right)^{1-y}\frac{a}{b}\leq \frac{B(b,y)}{B(a,y)}\leq\left(\frac{b}{a}\right)^{-y} .$$
In the case $y>1$, since $(y-1)>0$, one has to use instead
$$-(y-1)\int_0^1 t^{x-1}(1-t)^{y-1} dt\leq I\leq -(y-1)\int_0^1 t^{x}(1-t)^{y-1} dt.$$
The proof is complete.
\end{proof}
\blem\label{lem:compare1}
The following assertions are satisfied.
\begin{itemize}
\item[(a)] For $0<a<b$ and $0<y<\infty$, we have 
\Be\label{eq:compare11}
\left(\frac{a}{b}\right)^{1+y}\left(\frac{b+y}{a+y}\right)<\frac{b^{b-1}(a+y)^{a+y-1}}{a^{a-1}(b+y)^{b+y-1}}<\left(\frac{b+y}{a+y}\right)^{1-y}\left(\frac{a}{b}\right).
\Ee
%and
%\Be\label{eq:compare12}\left(\frac{a}{b}\right)^y\leq \frac{b^{b}(a+y)^{a+y}}{a^{a}(b+y)^{b+y}}.\Ee

\item[(b)] For $0<y<1$, and $\frac{1}{y}-y<a<b$, we have \Be\label{eq:compare13}\left(\frac{a}{b}\right)^y\leq \frac{b+y}{a+y}\left(\frac{a}{b}\right)\left(\frac{a+y+1}{b+y+1}\right)^y.\Ee
\item[(c)] We have that for $0<a<b<\frac{y}{y-2}$, and $y\geq 2$,
\Be\label{eq:compare14}\left(\frac{b+y}{a+y}\right)^{1-y}\left(\frac{a}{b}\right)< \frac{b^{b}(a+y)^{a+y}}{a^{a}(b+y)^{b+y}}.\Ee
\item[(d)] For $\frac{y^2-y}{2-y}<a<b$ with $1<y<2$,
\Be\label{eq:compare15}\left(\frac{a}{b}\right)^y\geq \frac{b^{b-1}(a+y)^{a+y-1}}{a^{a-1}(b+y)^{b+y-1}}.\Ee
\end{itemize}
\elem
\begin{remark}\label{rem:main1}
\begin{itemize}
\item Assertion (a) above implies that (\ref{eq:main11}) improves (\ref{eq:Fromlowerupper}) when $0<y \leq 1$.
\medskip

\item Assertion (b) shows that (\ref{eq:main11}) also improves (\ref{eq:twoside1}) at least for $0<y<1$ and $\frac{1}{y}-y<a<b$.
\item From assertion (c), we deduce that
the upper bound in (\ref{eq:main12}) improves the one in (\ref{eq:Fromlowerupper}) for $0<a<b<\frac{y}{y-2}$ with $y\geq 2$.

\item Assertion (d) implies that for $\frac{y^2-y}{2-y}<a<b$ with $1<y<2$,
the lower bound in (\ref{eq:main12}) is better than the one in (\ref{eq:Fromlowerupper}).

\item It is obvious that $$\left(\frac{b+y}{a+y}\right)^{1-y}\left(\frac{a}{b}\right)\leq \frac{b+y}{a+y}\left(\frac{a}{b}\right)\left(\frac{a+y+1}{b+y+1}\right)^y.$$  
Thus the upper bound in (\ref{eq:main12}) improves the one in (\ref{eq:twoside1}).
\end{itemize}
\end{remark}
\begin{proof}[Proof of Lemma \ref{lem:compare1}]
The right inequality in (\ref{eq:compare11}) follows from the fact that the function $x\mapsto x\ln\left(1+\frac{y}{x}\right)$ is increasing.
Observing that the function $x\mapsto (x+y)\left[\ln(x+y)-\ln x\right]$ is decreasing, one deduces the left inequality. Inequality (\ref{eq:compare13}) is obtained using the same argument.

Let us prove assertions (c) and (d). We start by observing that the inequality (\ref{eq:compare14}) is equivalent to $$1\leq \frac{A(b)}{A(a)}\quad\text{where}\quad A(x)=\frac{x^{x+1}}{(x+y)^{x+1}}.$$
Put $$f(x)=\ln A(x)=(x+1)\left[\ln x-\ln(x+y)\right].$$
Then $$f'(x)=\ln x-\ln(x+y)+\frac{1}{x}+\frac{y-1}{x+y}$$
and $$f''(x)=\frac{y(xy-2x-y)}{x^2(x+y)^2}.$$
For $y>2$, one has that $f''(x)<0$ for all $0<x<\frac{y}{y-2}$. Thus $f'$ is decreasing on $(0,\frac{y}{y-2})$ and as $f'\left(\frac{y}{y-2}\right)\geq 0$, we conclude that $f'(x)>0$. Thus $f(x)$ is increasing. Hence $f(b)>f(a)$ or equivalently, $\frac{A(b)}{A(a)}>1$. 
\medskip

As above, we have that the inequality (\ref{eq:compare15}) is equivalent to $$1\geq \frac{A(b)}{A(a)}\quad\text{where}\quad A(x)=\left(\frac{x}{x+y}\right)^{x+y-1}.$$  
Again, take
$$f(x)=\ln A(x)=(x+y-1)\left(\ln x-\ln(x+y)\right).$$
Then $$f'(x)=\ln x-\ln(x+y)+\frac{y-1}{x}+\frac{1}{x+y}$$
and $$f''(x)=\frac{y(y-yx-y^2+2x)}{x^2(x+y)^2}.$$
Then as $f''(x)>0$ for $x>\frac{y^2-y}{2-y}$ for $1<y<2$. Also as $\displaystyle\lim_{x\mapsto\infty}f'(x)=0$,  we conclude that $f'(x)<0$, and so $f(x)$ is decreasing. Thus $\frac{A(b)}{A(a)}<1$. The proof is complete.
\end{proof}
We have the following.
\btheo\label{thm:main2}
Let $1<a<b$. Then 
    \Be\label{eq:main22}
 e^{2(a-b)}\left(\frac{a-1}{b-1}\right)\leq \frac{B(b,b)}{B(a,a)}\leq e^{a-b}\left(\frac{a}{b}\right)^{\frac{1}{2}}.
 \Ee
\etheo
\begin{proof}
 In (\ref{eq:integparts3}), make $y=x$, to obtain
 \Be\label{eq:integparts31}
 \int_0^1 t^{x-1}(1-t)^{x-1} dt+x\int_0^1 t^{x-1}(1-t)^{x-1}\left(\ln t\right) dt=(x-1)\int_0^1 t^{x}(1-t)^{x-2}\left(\ln t\right) dt.
 \Ee
 We have by the product rule that
 \Beas
 \partial_xB(x,x) &=& \int_0^1 t^{x-1}(1-t)^{x-1}\left[\ln t+\ln(1-t)\right] dt.
 \Eeas
 That is
 \Be\label{eq:diffBxx}
 \int_0^1 t^{x-1}(1-t)^{x-1}\left(\ln t\right) dt=\partial_xB(x,x)-\int_0^1 t^{x-1}(1-t)^{x-1}\left(\ln(1- t)\right)dt.
 \Ee
 Putting together (\ref{eq:integparts3}) and (\ref{eq:diffBxx}), we obtain that
 \Be\label{eq:diffBxxx}
 B(x,x)+x\partial_xB(x,x)=x\int_0^1 t^{x-1}(1-t)^{x-1}\left(\ln(1- t)\right)dt+(x-1)\int_0^1 t^{x}(1-t)^{x-2}\left(\ln(t)\right)dt.
 \Ee
 Let us put
 $$I:=x\int_0^1 t^{x-1}(1-t)^{x-1}\left(\ln(1- t)\right)dt+(x-1)\int_0^1 t^{x}(1-t)^{x-2}\left(\ln(t)\right)dt.$$
 Using again the inequalities $\frac{(t-1)}{t}\leq \ln t\leq t-1$ for all $0<t\leq 1$, we obtain for all $x>1$,
$$-\int_0^1t^x(1-t)^{x-2}dt\leq \int_0^1t^{x-1}(1-t)^{x-1}\left(\ln(1-t)\right)dt\leq -\int_0^1t^x(1-t)^{x-1}dt$$
and 
 $$-\int_0^1t^{x-1}(1-t)^{x-1}dt\leq \int_0^1t^{x}(1-t)^{x-2}\left(\ln t\right)dt\leq -\int_0^1t^x(1-t)^{x-1}dt.$$
 It follows that 
 $$-xB(x+1,x-1)-(x-1)B(x,x)\leq I\leq -xB(x+1,x)-(x-1)B(x+1,x)$$
 or equivalently,
 \Be\label{eq:estimI}
 \left(-\frac{x^2}{x-1}-x+1\right)B(x,x)\leq I\leq \left(-x+\frac{1}{2}\right)B(x,x).
 \Ee
 From (\ref{eq:estimI}) and (\ref{eq:diffBxxx}), we deduce that 
 $$\left(-\frac{x^2}{x-1}-x+1\right)B(x,x)\leq  B(x,x)+x\partial_xB(x,x)\leq \left(-x+\frac{1}{2}\right)B(x,x).$$
 This is equivalent to
 
 \Be\label{eq:diffquotient2}
-\frac{1}{x-1}-2\leq \frac{\partial_xB(x,x)}{B(x,x)}\leq -1-\frac{1}{2x}.
\Ee
Integrating from $x=a$ to $x=b$, we obtain
$$-2(b-a)-\ln\left(\frac{b-1}{a-1}\right)\leq \ln\left(\frac{B(b,b)}{B(a,a)}\right)\leq -(b-a)-\frac{1}{2}\ln\left(\frac{b}{a}\right).$$
This is equivalent to (\ref{eq:main22}). Inequalities (\ref{eq:main22}). The proof is complete.
\end{proof}
We have seen in Remark \ref{rem:main1} that the upper bound in (\ref{eq:main12}) improves the one in (\ref{eq:Fromlowerupper}) for $0<a<b<\frac{y}{y-2}$ with $y\geq 2$, and that the lower bound in (\ref{eq:main12}) is better than the one in (\ref{eq:Fromlowerupper}) for $\frac{y^2-y}{2-y}<a<b<\infty$ with $1<y<2$. To improve these ranges of the parameters $a$ and $b$, we prove the following.
\btheo\label{thm:main3}
Let $0<a<b$. If $y\geq 1$, then
\Be\label{eq:main3}
\frac{b^{b-1}(a+y-1)^{a+y-1}}{a^{a-1}(b+y-1)^{b+y-1}}\leq\frac{B(b,y)}{B(a,y)}\leq \frac{a(b+1)^{b+1}(a+y)^{a+y}}{b(a+1)^{a+1}(b+y)^{b+y}}.
\Ee
\etheo
\begin{proof}
 Let $y>1$. Using (\ref{eq:integparts2}), we obtain
 $$B(x,y)+x\partial_xB(x,y)=(y-1)\left(\int_0^1\frac{t^{x-1}(1-t)^{y-2}}{B(x,y-1)}(t\ln t)dt\right)B(x,y-1).$$
 We note that the function $t\mapsto \Phi(t)=t\ln t$ is convex. We also note that $$d\mu(t)=\frac{t^{x-1}(1-t)^{y-2}}{B(x,y-1)}dt$$ is a probability measure.
 \medskip

 Using the Jensen's inequality, we obtain
 \Beas
 \int_0^1\frac{t^{x-1}(1-t)^{y-2}}{B(x,y-1)}(t\ln t)dt &=& \int_0^1\Phi(t)d\mu(t)\\ &\geq& \Phi\left(\int_0^1td\mu(t)\right)\\ &=& \Phi\left(\int_0^1\frac{t^x(1-t)^{y-2}}{B(x,y-1)}dt\right)\\ &=& \Phi\left(\frac{B(x+1,y-1)}{B(x,y-1)}\right)\\ &=& \Phi\left(\frac{x}{x+y-1}\right)\\ &=& \frac{x}{x+y-1}\ln\left(\frac{x}{x+y-1}\right).
 \Eeas
 Thus
 \Beas
 B(x,y)+x\partial_xB(x,y) &\geq& (y-1)\frac{x}{x+y-1}\ln\left(\frac{x}{x+y-1}\right)B(x,y-1)\\ &=& x\ln\left(\frac{x}{x+y-1}\right)B(x,y).
 \Eeas
 Consequently,
 \Be\label{eq:diffmain31}\frac{\frac{\partial B(x,y)}{\partial x}}{B(x,y)}\geq \ln\left(\frac{x}{x+y-1}\right)-\frac 1x.\Ee
 Integrating both sides from $x=a$ to $x=b$, we obtain
 $$\ln\left(\frac{B(b,y)}{B(a,y)}\right)\geq b\ln b-(b+y-1)\ln(b+y-1)-a\ln a+(a+y-1)\ln(a+y-1)-\ln b+\ln a.$$
 Thus $$\frac{B(b,y)}{B(a,y)}\geq \frac{b^{b-1}}{a^{a-1}}\frac{(a+y-1)^{a+y-1}}{(b+y-1)^{b+y-1}}.$$
 
 Let us prove the upper bound. We use the fact that the function $\Psi(t)=\ln t$ is concave, and $$d\nu(t)=\frac{t^{x}(1-t)^{y-2}}{B(x+1,y-1)}dt$$ is a probability measure.  Using Jensen's inequality, we obtain
 \Beas
 \int_0^1\frac{t^{x}(1-t)^{y-2}}{B(x+1,y-1)}(\ln t)dt &=& \int_0^1\Psi(t)d\nu(t)\\ &\leq& \Psi\left(\int_0^1td\nu(t)\right)\\ &=& \Psi\left(\int_0^1\frac{t^{x+1}(1-t)^{y-2}}{B(x+1,y-1)}dt\right)\\ &=& \Psi\left(\frac{B(x+2,y-1)}{B(x+1,y-1)}\right)\\ &=& \Psi\left(\frac{x+1}{x+y}\right)\\ &=& \ln\left(\frac{x+1}{x+y}\right).
 \Eeas
 Thus
 \Beas
 B(x,y)+x\partial_xB(x,y) &\leq& (y-1)\ln\left(\frac{x+1}{x+y}\right)B(x+1,y-1)\\ &=& x\ln\left(\frac{x+1}{x+y}\right)B(x,y).
 \Eeas
 It follows that
 \Be\label{eq:diffmain32}\frac{\frac{\partial B(x,y)}{\partial x}}{B(x,y)}\leq \ln\left(\frac{x+1}{x+y}\right)-\frac 1x.\Ee
 Integrating both sides from $x=a$ to $x=b$, we obtain
 $$\ln\left(\frac{B(b,y)}{B(a,y)}\right)\leq -\ln b+\ln a+(b+1)\ln(b+1)-(a+1)\ln(a+1)-(b+y)\ln(b+y)+(a+y)\ln(a+y).$$
 Thus $$\frac{B(b,y)}{B(a,y)}\leq \frac{a(b+1)^{b+1}(a+y)^{a+y}}{b(a+1)^{a+1}(b+y)^{b+y}}.$$
 The proof is complete.
\end{proof}
\begin{remark}\label{rem:main3}
The above result improves  both (\ref{eq:Fromlowerupper}) and (\ref{eq:twoside1}) for this range of the parameters. Also, Theorem \ref{thm:main3} improves Theorem \ref{thm:main1}. Indeed, we have the following inequalities.
\blem\label{lem:compare2}
For $0<a<b$ and $y\geq 1$, the following inequalities hold.\Be\label{eq:main3vsfromUG1}\frac{a(b+1)^{b+1}(a+y)^{a+y}}{b(a+1)^{a+1}(b+y)^{b+y}}\leq \frac{a(b+y)}{b(a+y)}\left(\frac{a+y+1}{b+y+1}\right)^y;\Ee  

\Be\label{eq:main3vsfromUG2}\frac{b^{b-1}(a+y-1)^{a+y-1}}{a^{a-1}(b+y-1)^{b+y-1}}\geq\frac{b^{b-1}(a+y)^{a+y-1}}{a^{a-1}(b+y)^{b+y-1}};\Ee

\Be\label{eq:main3vsfromUG3}
\frac{a(b+1)^{b+1}(a+y)^{a+y}}{b(a+1)^{a+1}(b+y)^{b+y}}\leq\frac{b^b(a+y)^{a+y}}{a^a(b+y)^{b+y}};
\Ee
\Be\label{eq:main3vsfromUG4}\left(\frac{a}{b}\right)^{y}\leq \frac{b^{b-1}(a+y-1)^{a+y-1}}{a^{a-1}(b+y-1)^{b+y-1}},\Ee
and
\Be\label{eq:main3vsfromUG5}\frac{a(b+1)^{b+1}(a+y)^{a+y}}{b(a+1)^{a+1}(b+y)^{b+y}} \leq\left(\frac{b+y}{a+y}\right)^{1-y}\left(\frac{a}{b}\right).\Ee
\elem
\begin{proof}
Note that (\ref{eq:main3vsfromUG1}) is equivalent to $$A(a,b,y)\leq B(a,b,y)$$ where
$A(a,b,y)=\frac{(b+1)^{b+1}(a+y)^{a+y+1}}{(a+1)^{a+1}(b+y)^{b+y+1}}$ and $B(a,b,y)=\left(\frac{a+y+1}{b+y+1}\right)^y$. Now we have
$$\ln A(a,b,y)-\ln B(a,b,y)=f(b)-f(a)$$
where $f(x)=(x+1)\ln(x+1)-(x+y+1)\ln(x+y)+y\ln(x+y+1)$. It suffices then to see that $f(x)$ is decreasing to conclude. The proof of (\ref{eq:main3vsfromUG2}), (\ref{eq:main3vsfromUG3}), (\ref{eq:main3vsfromUG4}) and (\ref{eq:main3vsfromUG5}) follows in the same way.
\end{proof}
\end{remark}
The above result does not take into account the case $0<y\leq 1$. Here is how to obtain a related estimate in this case.
\btheo\label{thm:main4}
Let $0<a<b$. The following hold.
\begin{itemize}
    \item[(1)] If $0<y\leq 1$, then
    \Be\label{eq:main41}
    \frac{B(b,y)}{B(a,y)}\leq \frac{a}{b}\left(\frac{b^b(a+y)^{a+y}}{a^a(b+y)^{b+y}}\right)^{1-\frac{1}{y}}.
    \Ee
    \item[(2)] If $y>1$, then 
    \Be\label{eq:main42}
 \frac{a}{b}\left(\frac{b^b(a+y)^{a+y}}{a^a(b+y)^{b+y}}\right)^{1-\frac{1}{y}} \leq \frac{B(b,y)}{B(a,y)}.
    \Ee
\end{itemize}
\etheo
\begin{proof}
 We proceed as above. We note that the function $\Psi(t)=\frac{t\ln t}{1-t}$ is convex, and we recall that $$d\nu(t)=\frac{t^{x-1}(1-t)^{y-1}}{B(x,y)}dt$$ is a probability measure.  Using again Jensen's inequality, we obtain
 \Beas
 \int_0^1\frac{t^{x-1}(1-t)^{y-1}}{B(x,y)}\left(\frac{t\ln t}{1-t}\right)dt &=& \int_0^1\Psi(t)d\nu(t)\\ &\geq& \Psi\left(\int_0^1td\nu(t)\right)\\ &=& \Psi\left(\int_0^1\frac{t^{x}(1-t)^{y-1}}{B(x,y)}dt\right)\\ &=& \Psi\left(\frac{B(x+1,y)}{B(x,y)}\right)\\ &=& \Psi\left(\frac{x}{x+y}\right)\\ &=& \frac{x}{y}\ln\left(\frac{x}{x+y}\right).
 \Eeas
 It follows from (\ref{eq:integparts2}) and the above that in the case $0<y\leq 1$,
 \Beas
 B(x,y)+x\partial_xB(x,y) &\leq& \frac{(y-1)x}{y}\ln\left(\frac{x}{x+y}\right)B(x,y).
 \Eeas
 Thus
 \Be\label{eq:diffmain41}\frac{\frac{\partial B(x,y)}{\partial x}}{B(x,y)}\leq \left(1-\frac{1}{y}\right)\ln\left(\frac{x}{x+y}\right)-\frac 1x.\Ee
 Integrating both sides from $x=a$ to $x=b$, we obtain
 $$\ln\left(\frac{B(b,y)}{B(a,y)}\right)\leq -\ln b+\ln a+\left(1-\frac{1}{y}\right)\left(b\ln b-a\ln a-(b+y)\ln(b+y)+(a+y)\ln(a+y)\right).$$
 Thus $$\frac{B(b,y)}{B(a,y)}\leq \frac{a}{b}\left(\frac{b^b(a+y)^{a+y}}{a^a(b+y)^{b+y}}\right)^{1-\frac{1}{y}}.$$ 
 The case $y>1$ follows similarly, but with the reverse inequality due to $(y-1)>0$.
\end{proof}
\begin{remark}\label{rem:main4}
We have that (\ref{eq:main41}) improves the upper bound in (\ref{eq:main11}) for $0<y<1$, and that (\ref{eq:main42}) improves the lower bound in (\ref{eq:main3}) for all  $y\geq 1$ and $0<a<b$. This is a consequence of the following result.
\blem\label{lem:compare3}
Let   $0<a<b$, then the following inequalities hold. For $0<y<1$, we have 
\Be\label{eq:main4vsfromUG2}\frac{a}{b}\left(\frac{b^b(a+y)^{a+y}}{a^a(b+y)^{b+y}}\right)^{1-\frac{1}{y}}\leq \left(\frac{a}{b}\right)^y,\Ee
and for $y \geq 1$, we have
\Be\label{eq:main4vsfromUG3}\frac{a}{b}\left(\frac{b^b(a+y)^{a+y}}{a^a(b+y)^{b+y}}\right)^{1-\frac{1}{y}}\geq \frac{b^{b-1}(a+y-1)^{a+y-1}}{a^{a-1}(b+y-1)^{b+y-1}}.\Ee
\elem
\begin{proof}
 Inequality (\ref{eq:main4vsfromUG2}) can be proved as in Lemma \ref{lem:compare2}.

\noindent
Proving inequality \eqref{eq:main4vsfromUG3} is equivalent to showing that $\frac{A(a)}{A(b)} \geq 1$,
where the function $A$ is defined by
\[
A(x) = \frac{x^{\frac{x}{y}}(x+y)^{(x+y)\left(1-\frac{1}{y}\right)}}{(x+y-1)^{x+y-1}}.
\]
Let
\[
f(x) = \ln A(x) = \frac{x}{y}\ln x + \left(1-\frac{1}{y}\right)(x+y)\ln(x+y) - (x+y-1)\ln(x+y-1).
\]
It suffices to prove that $f$ is decreasing on $(0,\infty)$. We have
\[
f'(x) = \frac{1}{y}\ln\!\left(\frac{x}{x+y}\right) + \ln\!\left(\frac{x+y}{x+y-1}\right),
\]
and
\[
f''(x) = \frac{y-1}{x(x+y)(x+y-1)}.
\]
Since $y \geq 1$, it follows that $f''(x) \geq 0$ for all $x > 0$, hence $f'$ is increasing on $(0,\infty)$. Moreover,
\[
\lim_{x \to \infty} f'(x) = 0.
\]
Therefore, $f'(x) \leq 0$ for all $x > 0$, which shows that $f$ is decreasing on $(0,\infty)$.
\end{proof}
\end{remark}

We have the following result.
\btheo\label{thm:main6}
Assume that $0<a<b<\infty$, $y>0$. Let 
$$c=1-\frac{y}{b-a}\ln\left(\frac{b+y+1}{a+y+1}\right),\quad \alpha=\ln c-1,\quad \beta=\frac 1c.$$
Then
\Be\label{eq:main6}
\frac{B(b,y)}{B(a,y)}\leq e^{(\alpha+\beta)(b-a)}\left(\frac{a+y+1}{b+y+1}\right)^{y\beta}\frac{a}{b}\left(\frac{b+y}{a+y}\right)^y.
\Ee
\etheo
\begin{proof}
 We take $k=l$ and $h=0$ in (\ref{eq:integparts1}) to get  
 \Be\label{eq:integparts6}
 B(x,y+1)+x\partial_xB(x,y+1)=y\int_0^1 t^{x}(1-t)^{y-1}(\ln t) dt.
 \Ee
 We now apply an idea from \cite{From}. That is we note that as $t\mapsto \ln t$ is concave, for any $0<t<1$ and any $c>0$,
 $$\ln t\leq \ln c-1+\frac{t}{c}=\alpha+\beta t.$$
 It follows that
 \Beas
 \int_0^1 t^{x}(1-t)^{y-1}(\ln t) dt &\leq& \alpha B(x+1,y)+\beta B(x+2,y)\\ &=& \left(\alpha\frac{x}{x+y}+\beta\frac{x(x+1)}{(x+y)(x+y+1)}\right)B(x,y).
 \Eeas
 Taking this into (\ref{eq:integparts6}), we obtain
 $$B(x,y+1)+x\partial_xB(x,y+1)\leq y\left(\alpha\frac{x}{x+y}+\beta\frac{x(x+1)}{(x+y)(x+y+1)}\right)B(x,y)$$
 or equivalently,
 \Be\label{eq:diffmain6}
 \frac{\frac{\partial B(x,y)}{\partial x}}{B(x,y)}\leq \alpha+\beta-\frac{y\beta}{x+y+1}-\frac{1}{x}+\frac{y}{x+y}.
 \Ee
 Integrating both sides from $x=a$ to $x=b$, we obtain
 \Be\label{eq:Ineqlog6}\ln\left(\frac{B(b,y)}{B(a,y)}\right)\leq (\alpha+\beta)(b-a)-y\beta\ln\left(\frac{b+y+1}{a+y+1}\right)-\ln\left(\frac{b}{a}\right)+y\ln\left(\frac{b+y}{a+y}\right).
 \Ee
 One easily checks that the value of $c>0$ minimizing (\ref{eq:Ineqlog6}) is given as 
 $$c=1-\frac{y}{b-a}\ln\left(\frac{b+y+1}{a+y+1}\right).$$
 Note that $0<c<1$. Taking this into (\ref{eq:Ineqlog6}) and then taking the exponential of both members, we obtain (\ref{eq:main6}). The proof is complete.
\end{proof}
\begin{remark}\label{rem:main6}
 Let us compare the bounds in Theorem \ref{thm:From317} and Theorem \ref{thm:main6}.  
 \begin{figure}[H]
 \centering
 \includegraphics[width=0.9\textwidth]{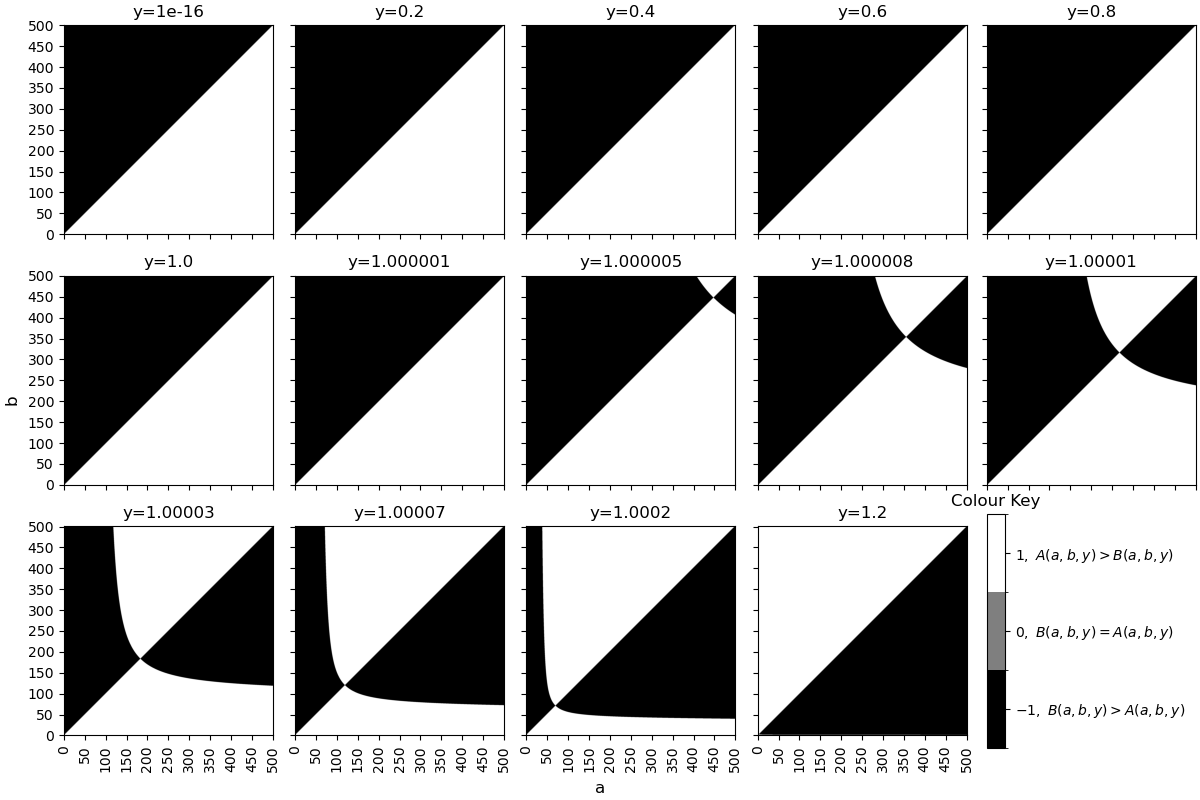}
 \caption{Plot showing the sign of $F(a,b,y)$ . The vertical axis corresponds to $b$ and the horizontal axis corresponds to $a$.  \label{figlarge}}
 \end{figure}
 %showing the parts of the plane where $A(a,b,y)>B(a,b,y)$ in white and $A(a,b,y)<B(a,b,y)$ in black. $y$ is taken from $10^{-16}$ to 1.2
 We define the function
 \begin{equation}
 	F(a,b,y) = sign(A(a,b,y)-B(a,b,y))=\begin{cases}
 		1 & \text{if }\ A(a,b,y)-B(a,b,y) > 0\\
 		0 & \text{if }\ A(a,b,y)-B(a,b,y) = 0\\
 		-1 & \text{if }\ A(a,b,y)-B(a,b,y) < 0.
 	\end{cases}
 \end{equation} 
 where $$A\left( a,b,y\right)= e^{(\alpha+\beta)(b-a)}\left(\frac{a+y+1}{b+y+1}\right)^{y\beta}\frac{a}{b}\left(\frac{b+y}{a+y}\right)^y$$ from Theorem \ref{thm:main6} and 
 $$B\left( a,b,y\right)=  e^{(\alpha+\beta)(b-a)}\left(\frac{a+y}{b+y}\right)^{y\beta} $$ from Theorem \ref{thm:From317}. We can see from figure \ref{figlarge} that for $0<y\le1,$ $A(a,b,y)<B(a,b,y)$  if $b>a$ which is the portion of the plane above the diagonal. From $1.000001<y<1.2,$ this initial behaviour slowly flips over. As $y\rightarrow 1.2$ the area in the  $a-b$ plane with $b>a$ where $A(a,b,y)<B(a,b,y)$ keeps getting smaller. When $y>1.2,$ when $b>a,$ we have that $A(a,b,y)>B(a,b,y).$ 
 
 %\clearpage
\end{remark}
Let us prove the following lower estimate.
\btheo\label{thm:main7}
Assume that $0<a<b<\infty$, $y>0$. Let 
$$c=\frac{b-a}{b-a+y\ln\left(\frac{b}{a}\right)},\quad \alpha=\ln c+1,\quad \beta=-c.$$
Then
\Be\label{eq:main7}
\frac{B(b,y)}{B(a,y)}\geq e^{(\alpha+\beta)(b-a)}\left(\frac{b}{a}\right)^{y\beta-1}\left(\frac{b+y}{a+y}\right).
\Ee
\etheo
\begin{proof}
We rewrite (\ref{eq:integparts6}) as   
\Be\label{eq:integparts7}
 B(x,y+1)+x\partial_xB(x,y+1)=y\int_0^1 t^{x-1}(1-t)^{y-1}(t\ln t) dt.
 \Ee
 We note that as $t\mapsto t\ln t$ is convex, for any $0<t<1$ and any $c>0$,
 $$t\ln t\geq (\ln c+1)t-c=\alpha t+\beta .$$
 It follows that
 \Beas
 \int_0^1 t^{x-1}(1-t)^{y-1}(t\ln t) dt &\geq& \alpha B(x+1,y)+\beta B(x,y)\\ &=& \left(\alpha\frac{x}{x+y}+\beta\right)B(x,y).
 \Eeas
Taking this into (\ref{eq:integparts7}), we obtain
 $$B(x,y+1)+x\partial_xB(x,y+1)\geq y\left(\alpha\frac{x}{x+y}+\beta\right)B(x,y)$$
 or equivalently,
 \Be\label{eq:diffmain7}
 \frac{\frac{\partial B(x,y)}{\partial x}}{B(x,y)}\geq (\alpha+\beta)+\beta\frac{y}{x}-\frac{1}{x}+\frac{1}{x+y}.
 \Ee
 Integrating both sides from $x=a$ to $x=b$, we obtain
 \Be\label{eq:Ineqlog7}\ln\left(\frac{B(b,y)}{B(a,y)}\right)\geq (\alpha+\beta)(b-a)+y\beta\ln\left(\frac{b}{a}\right)-\ln\left(\frac{b}{a}\right)+\ln\left(\frac{b+y}{a+y}\right).
 \Ee
 We have that the value of $c>0$ maximizing (\ref{eq:Ineqlog7}) is given as 
 $$\frac{1}c=1+\frac{y}{b-a}\ln\left(\frac{b}{a}\right).$$
 Taking this into (\ref{eq:Ineqlog7}) and then taking the exponential of both members, we obtain (\ref{eq:main7}). The proof is complete.
\end{proof}
The following improves the lower bound in Theorem \ref{thm:main2}.
\btheo\label{thm:main22}
Let $0<a<b$. Then the following holds.
    \Be\label{eq:main221}
 e^{2(a-b)}\left(\frac{a}{b}\right)\leq \frac{B(b,b)}{B(a,a)}\leq e^{a-b}\left(\frac{a}{b}\right)\left(\frac{2b+1}{2a+1}\right)^{\frac{1}{2}}.
    \Ee
\etheo
\begin{proof}
  Take $k=l$ and $h=0$ in (\ref{eq:integparts1}), and then make $y=x$, to obtain
 \Be\label{eq:integparts01}
 \int_0^1 t^{x-1}(1-t)^{x} dt+x\int_0^1 t^{x-1}(1-t)^{x}\left(\ln t\right) dt=x\int_0^1 t^{x}(1-t)^{x-1}\left(\ln t\right) dt.
 \Ee
 We have by the product rule that
 \Beas
 \partial_xB(x,x+1) &=& \dfrac{d}{dx}\int_0^1 t^{x-1}(1-t)^{x} dt\\ &=&\int_0^1 t^{x-1}(1-t)^{x}\left[\ln t+\ln(1-t)\right] dt.
 \Eeas
 That is
 \Be\label{eq:diffBxx1}
 \int_0^1 t^{x-1}(1-t)^{x}\left(\ln t\right) dt=\partial_xB(x,x+1)-\int_0^1 t^{x-1}(1-t)^{x}\left(\ln(1- t)\right)dt.
 \Ee
 Putting together (\ref{eq:integparts01}) and (\ref{eq:diffBxx1}), we obtain that
 \Be\label{eq:diffBxxx1}
 B(x,x+1)+x\partial_xB(x,x+1)=x\int_0^1 t^{x-1}(1-t)^{x}\left(\ln(1- t)\right)dt+x\int_0^1 t^{x}(1-t)^{x-1}\left(\ln(t)\right)dt.
 \Ee
 Let us put
 $$I:=x\int_0^1 t^{x-1}(1-t)^{x}\left(\ln(1- t)\right)dt+x\int_0^1 t^{x}(1-t)^{x-1}\left(\ln(t)\right)dt.$$
 As above, we use that $\frac{(t-1)}{t}\leq \ln t\leq t-1$ for all $0<t\leq 1$, to obtain for all $x>0$,
$$-\int_0^1t^x(1-t)^{x-1}dt\leq \int_0^1t^{x-1}(1-t)^{x}\left(\ln(1-t)\right)dt\leq -\int_0^1t^x(1-t)^{x}dt$$
and 
 $$-\int_0^1t^{x-1}(1-t)^{x}dt\leq \int_0^1t^{x}(1-t)^{x-1}\left(\ln t\right)dt\leq -\int_0^1t^x(1-t)^{x}dt.$$
 It follows that 
 \Be\label{eq:estimI1}
 -2xB(x,x+1)\leq I\leq -2xB(x+1,x+1)=-\frac{2x^2}{2x+1}B(x,x+1).
 \Ee
 From (\ref{eq:estimI1}) and (\ref{eq:diffBxxx1}), we deduce that 
 $$-2xB(x,x+1)\leq  B(x,x+1)+x\partial_xB(x,x+1)\leq -\frac{2x^2}{2x+1}B(x,x+1).$$
 This is equivalent to
 
 \Be\label{eq:diffquotient200}
\left(-2-\frac{1}{x}\right)\leq \frac{\partial_xB(x,x)}{B(x,x)}\leq \left(-\frac{2x}{2x+1}-\frac{1}{x}\right)=-1+\frac{1}{2x+1}-\frac{1}{x}.
\Ee
Integrating from $x=a$ to $x=b$, we obtain
%$$-2(b-a)-\ln\left(\frac{b}{a}\right)\leq \ln\left(\frac{B(b,b+1)}{B(a,a+1)}\right)\leq -(b-a)+\frac{1}{2}\ln\left(\frac{2b+1}{2a+1}\right)-\ln\left(\frac{b}{a}\right)$$
%or equivalently,
$$-2(b-a)-\ln\left(\frac{b}{a}\right)\leq \ln\left(\frac{B(b,b)}{B(a,a)}\right)\leq -(b-a)+\frac{1}{2}\ln\left(\frac{2b+1}{2a+1}\right)-\ln\left(\frac{b}{a}\right).$$
This is equivalent to (\ref{eq:main221}). The proof is complete.
\end{proof}
\begin{remark}
For $1<a<b<\infty$, it is easy to see that $$\frac{a-1}{b-1}<\frac{a}{b}$$
while

$$\left(\frac{2b+1}{2a+1}\right)^{1/2}\left(\frac{a}{b}\right)<\left(\frac{a}{b}\right)^{1/2}.$$
Hence, for $1<a<b$, the lower bound in (\ref{eq:main221}) improves the one in (\ref{eq:main22}) while the upper bound in (\ref{eq:main22}) is better than the one in (\ref{eq:main221}).
\end{remark}
The following result gives an estimate of the differences of the beta function.

\btheo\label{thm:main5}
The following hold.
\begin{itemize}
    \item[(1)] If $1\leq y\leq 2$, and $1\leq a<b\leq 2$, then 
    \Be\label{eq:main51}
B(b+1,y)\geq B(a+1,y) +\frac{1}{8\ln\left(\frac{3}{2}\right)}\left(\frac{5}{9}\right)^{y-1}\left[\left(\frac{4}{9}\right)^{b-1}-\left(\frac{4}{9}\right)^{a-1}\right].
\Ee
\item[(2)] If $\left(y>2\quad\text{and}\quad 0<a<b\leq 1\right)$ or  $\left(0<y<1\quad\text{and}\quad 2\leq a<b\right)$, then
\Be\label{eq:main52}
B(b+1,y)\leq B(a+1,y) +\frac{1}{8\ln\left(\frac{3}{2}\right)}\left(\frac{5}{9}\right)^{y-1}\left[\left(\frac{4}{9}\right)^{b-1}-\left(\frac{4}{9}\right)^{a-1}\right].
\Ee
\end{itemize}

\etheo
\begin{proof}
It is easy to verify that the function $\Phi_{x,y}(t)=t^{x-1}(1-t)^{y-1}$ defined in $(0,1)$, is convex if $(x,y)\in (0,1]\times[2,\infty)\cup[2,\infty)\times(0,1]$. It is concave if $(x,y)\in[1,2]\times[1,2]$. 
\medskip
Now consider the function $$g(x,y)=\int_0^1 t^x(1-t)^{y-1}dt.$$ 
Then $$\frac{\partial g(x,y)}{\partial x}=-\frac{1}{4}\int_0^1 t^{x-1}(1-t)^{y-1}(-4t\ln t)dt.$$
We observe that $(-4t\ln t)dt$ is a probability measure on $(0,1)$. Hence for $1\leq y\leq 2$ and $1\leq a<b\leq 2$, as $t\mapsto \Phi(t)=t^{x-1}(1-t)^{y-1}$ is concave, we obtain
\Beas
\frac{\partial g(x,y)}{\partial x} &=& -\frac{1}{4}\int_0^1\Phi(t) (-4t\ln t)dt\\ &\geq& -\frac{1}{4}\Phi\left(\int_0^1-4t^2\ln t dt\right)\\ &=& -\frac{1}{4}\Phi\left(\frac{4}{9}\right)\\ &=& -\frac{1}{4}\left(\frac{4}{9}\right)^{x-1}\left(\frac{5}{9}\right)^{y-1}.
\Eeas
That is $$\frac{\partial g(x,y)}{\partial x}\geq -\frac{1}{4}\left(\frac{5}{9}\right)^{y-1}e^{(x-1)\ln\left(\frac{4}{9}\right)}.$$
Integrating from $x=a$ to $x=b$, we obtain
$$B(b+1,y)-B(a+1,y)\geq \frac{1}{4\ln\left(\frac{9}{4}\right)}\left[\left(\frac{4}{9}\right)^{b-1}-\left(\frac{4}{9}\right)^{a-1}\right]\left(\frac{5}{9}\right)^{y-1}.$$
That (\ref{eq:main51}) holds.
\medskip

The estimate (\ref{eq:main52}) follows in the same way using that $t\mapsto \Phi(t)=t^{x-1}(1-t)^{y-1}$ is convex for the given ranges of $x$ and $y$.
\end{proof}
Taking $b=a+1$ in (\ref{eq:main52}), we deduce the following.
\bcor\label{cor:main4}
The following hold.
%\begin{itemize}
%\item[(a)] 
If $0<y<1$ and $a\geq 2$, then
\Be\label{eq:main52cor}
B(a,y)\geq \frac{1}{8\ln\left(\frac{3}{2}\right)}\left(\frac{5}{9}\right)^{y}\left(\frac{4}{9}\right)^{a-1}\frac{(a+y)(a+y+1)}{ay}.
\Ee
%\item[(b)] If $1\leq y,b\leq 2$, then
%\Be\label{eq:main51cor}
%B(b,y)\leq \frac{1}{8\ln\left(\frac{3}{2}\right)}\left(\frac{5}%{9}\right)^{y}\left(\frac{4}%{9}\right)^{b-2}\left(1+\frac{b}{y}\right).
%\Ee
%\end{itemize}
\ecor
%\begin{remark}
%Observe that for $0<y<1$ and $x\geq 3$, 
%$$\frac{1}{8\ln\left(\frac{3}{2}\right)}\left(\frac{5}%{9}\right)^{y}\left(\frac{4}%{9}\right)^{x-1}\frac{(x+y)(x+y+1)}{xy}\leq \frac{1}{xy}.$$
%Thus (\ref{eq:main52cor}) improves the upper bound in %(\ref{eq:dragomir}) for this range.
%\end{remark}

\begin{table}[h]
    \centering
    \begin{tabular}{|l|c|c|c|c|c|c|c|c|c|}
\hline 
 \multicolumn{2}{|c|}{\multirow{2}{*}{Theorems}} & \multicolumn{2}{c|}{$\dfrac{B(2,1)}{B(1,1)}=0.5$} & \multicolumn{2}{c|}{$\dfrac{B(5,8)}{B(4,8)}=0.333$} & \multicolumn{2}{c|}{$\dfrac{B(3,0.5)}{B(2,0.5)}=0.8$} & \multicolumn{2}{c|}{$\dfrac{B(5,5)}{B(4,5)}=0.4444$}\tabularnewline
%\hline
\cline{3-10} 
\multicolumn{2}{|c|}{}  & LB & UB & LB & UB & LB & UB & LB & UB\tabularnewline
\hline 
\hline 
\multirow{2}{*}{Existing} & \ref{thm:From} & 0.4444 & 0.5926 & 0.3114 & 0.3593 & 0.7762 & 0.8316 & 0.4204 & 0.4729\tabularnewline
\cline{2-10}
 & \ref{thm:twoside1} & 0.375 & 0.5625 & 0.1454 & 0.4790 & 0.7621 & 0.8231 & 0.2913 & 0.5519\tabularnewline
\hline 
\multirow{3}{*}{Proposed} & 2.2 & 0.5 & 0.5 & 0.1677 & 0.4568 & 0.7888 & 0.8165 & 0.3277 & 0.5249\tabularnewline
\cline{2-10}
 & 2.6 & 0.5 & 0.5 & 0.3125 & 0.3516 &  &  & 0.4229 & 0.4627\tabularnewline
\cline{2-10}
 & 2.9 &  & 0.5 & 0.3267 &  &  & 0.8017 & 0.4395 & \tabularnewline
\hline 
\end{tabular}\vspace{0.1cm}
    \caption{Numerical results comparing bounds in Theorems \ref{thm:From}, \ref{thm:twoside1}, 2.2, 2.6, and 2.9}
    \label{tab:ratios}
\end{table}

\begin{table}[h]
    \centering
    \begin{tabular}{|l|l|c|c|c|c|c|c|c|}
\hline 
\multicolumn{2}{|c|}{{\multirow{2}{*}{Theorems}}} & {\small B(2, 0.4) } & {\small B(2, 0.99) } & {\small B(2.5, 0.5) } & {\small B(2.8, 0.7) } & {\small B(3, 0.25) } & {\small B(4, 0.99) } & {\small B(5, 0.16) }\tabularnewline
\cline{3-9}
 \multicolumn{2}{|c|}{} & {\small\textbf{1.7857 }} & {\small\textbf{0.5076 }} & {\small\textbf{1.1781 }} & {\small\textbf{0.6548 }} & {\small\textbf{2.8444 }} & {\small\textbf{0.2553 }} & {\small\textbf{4.5541 }}\tabularnewline
\hline 
%\multirow{7}{*}{\rotatebox{90}{Existing}} & {\small 1.1} & {\small 1.6667 } & {\small 0.5067 } & {\small 1.0667 } & {\small 0.6033 } & {\small 2.4762 } & {\small 0.2541 } & {\small 3.5833 }\tabularnewline
%\cline{2-9}
\multirow{6}{*}{\rotatebox{90}{Existing}} & {\small \ref{thm:17}} & {\small 1.0174 } & {\small 0.2262 } & {\small 0.6211 } & {\small 0.3099 } & {\small 1.7949 } & {\small 0.1049 } & {\small 3.1608 }\tabularnewline
\cline{2-9}
 & {\small \ref{thm:31}} & {\small 0.0080 } &  &  &  & {\small 0.0156 } &  & {\small 0.0103 }\tabularnewline
\cline{2-9}
 & {\small \ref{thm:32}} & {\small 0.0476 } &  &  &  & {\small 0.0975 } &  & {\small 0.0703 }\tabularnewline
\cline{2-9}
 & {\small \ref{thm:34}} & {\small 0.0132 } &  &  &  & {\small 0.0259 } &  & {\small 0.0175 }\tabularnewline
\cline{2-9}
 & {\small \ref{thm:311}} &  &  &  &  & {\small 2.8444 } & {\small 0.1918 } & {\small 3.9632 }\tabularnewline
\cline{2-9}
 & {\small \ref{thm:314}} & {\small 0.1666 } & {\small 0.1073 } & {\small 0.0825 } & {\small 0.0488 } & {\small 0.0567 } & {\small 0.0095 } & {\small 0.0065 }\tabularnewline
\hline 
 & {\small Cor. 2.18} & {\small 1.1048 } & {\small 0.4614 } & {\small 0.6536 } & {\small 0.3814 } & {\small 0.9682 } & {\small 0.1142 } & {\small 0.4350 }\tabularnewline
\hline 
\end{tabular}\vspace{0.1cm}
    \caption{Numerical results comparing the lower bounds in Theorems \ref{thm:17} to \ref{thm:314} with Corollary 2.18}
    \label{tab:ineq}
\end{table}

\section{Numerical study}
In this section, we conduct a numerical study across various parameter ranges for $(a,b,y)$ and $(x,y)$. 
\medskip

It can be seen from Table \ref{tab:ratios} that Theorem \ref{thm:main4} gives the best bounds for all the values tested. This is closely followed by theorem \ref{thm:main3}. Theorem \ref{thm:main1} has mixed results, sometimes performing better or worse than both Theorem \ref{thm:From} and Theorem \ref{thm:twoside1}. From Table \ref{tab:ineq} it can be seen that Corollary \ref{cor:main4} generally provides very good bounds. The only ones that provide better bounds are Theorem \ref{thm:17} and Theorem \ref{thm:311}. However, Theorem \ref{thm:311} is only valid for $x\ge 3$ so Corollary \ref{cor:main4} gives bounds in the interval $2\le x <3$ where Theorem \ref{thm:311} does not. Theorem \ref{thm:17} only provides better bounds than Corollary \ref{cor:main4} in some circumstances.

%Table 1 compares the lower and upper bounds derived from our proposed theorems regarding the ratio of beta functions against existing results in the literature. The data demonstrates that, in specific cases, our results provide more accurate approximations for the ratio $\frac{B(b,y)}{B(a,y)}$ than those found in \cite{From}. Furthermore, Table 2 illustrates that our proposed upper bound for $B(x,y)$ outperforms the existing bounds established in \cite{From, Grenie} for certain values of $x$ and $y$. 

\section{Further remarks} 
In this section, we discuss some results that are implicit in the previous sections.
\medskip

Let us recall that the derivative of the logarithm of the  Euler's gamma function is called the digamma or psi function and denoted $\Psi(x)$. We note that to obtain our inequalities, we have first proved inequalities for $\partial_x\ln(B(x,y))=\frac{\partial_xB(x,y)}{B(x,y)}$. As we know that $$\partial_x\ln(B(x,y))=\Psi(x)-\Psi(x+y),$$ our method then also provides inequalities for the difference $\Psi(x+y)-\Psi(x)$. In the literature, estimates of $\Psi(x+y)-\Psi(x)$ are usually implicit, as they follow by integrating $\Psi'(t)$ from $t=x$ to $t=x+y$. For some estimates of $\Psi'(t)$ we refer the interested reader to \cite{Alzer3,Sun}.

Let us now state some estimates of $\Psi(x+y)-\Psi(x)$. From (\ref{eq:diffquotient1}) and (\ref{eq:diffmain41}), we deduce the following.
\bprop\label{prop:diffpsi1}
For $x>0$ and $0<y\leq 1$,
\Be\label{eq:diffpsi1}
\frac{1}{x}+\left(1-\frac{1}{y}\right)\ln\left(1+\frac{y}{x}\right)\leq \Psi(x+y)-\Psi(x)\leq \frac{y(x+1)}{x(x+y)}.
\Ee
\eprop
%It follows from the above and the First Principle of Calculus that the following holds.
%\bcor
%For $0<x<\infty$,
%\Be\label{eq:derivpsi}
%\frac{1}{x}\leq \Psi'(x)\leq \frac{1}{x}+\frac{1}%{x^2}.
%\Ee
%\ecor
From (\ref{eq:diffmain31}), (\ref{eq:diffmain32}) and the inequality obtained in place of (\ref{eq:diffmain41}) when $y>1$, we obtain the following.
\bprop\label{prop:diffpsi2}
For $x>0$ and $y>1$,
\Be\label{eq:diffpsi2}
\frac{1}{x}+\ln\left(\frac{x+y}{x+1}\right)\leq \Psi(x+y)-\Psi(x)\leq \frac{1}{x}+\left(1-\frac{1}{y}\right)\ln\left(1+\frac{y}{x}\right).
\Ee
\eprop
We now focus on the case $x=y$. In this direction, we note that in \cite{Alzer3}, H. Alzer observed that \Be\label{eq:Alzerpsi} \ln(x)-\frac{1}{x}\leq \Psi(x)\leq \ln(x)-\frac{1}{2x},\quad x>0.\Ee
BN. Guo and F. Qi proved in \cite[Theorem 1.1]{Guo} that for all $x>0$,
\Be\label{eq:Guo}\ln\left(x+\frac 12\right)-\frac{1}{x}\leq \Psi(x)\leq \ln\left(x+e^{-\gamma}\right)-\frac{1}{x}\Ee where $\gamma$ is the Euler-Mascheroni constant.

It is easy to check that the right hand side of (\ref{eq:Guo}) is less than the right hand side of (\ref{eq:Alzerpsi}). 

In \cite[Theorem 2.1]{Qiu}, the authors obtained that for $x>1$,
\Be\label{eq:Qiu}\frac{1}{2x}+\frac{1.5-\ln 4}{x^2}<\Psi\left(x+\frac 12\right)-\Psi(x)<\frac{1}{2x}+\frac{1}{8x^2}.\Ee
%For an upper bound for $\Psi\left(x+\frac 12\right)$, we have the following from \cite[Corollary]{Sun}
%\Be\label{eq:Sun} \Psi'\left(x+\frac 12\right)<\frac{1}%{3}\frac{4,800x^4+1,280x^2+147}{x(40x^2+7)^2}.\Ee
We show here how simpler bounds for $\Psi\left(x+\frac 12\right)-\Psi(x)$ can be derived from our proofs.
\medskip

Using inequalities (\ref{eq:diffquotient200}) and the identity $$\partial_x\ln(B(x,x))=\Psi(x)-\Psi\left(x+\frac 12\right)-\ln 4,$$ we obtain the following.
\bprop\label{prop:psi}
For $x>0$
    \Be\label{eq:psiany}
    1+\frac{1}{x}-\frac{1}{2x+1}-\ln 4 \leq \Psi\left(x+\frac 12\right)-\Psi(x)\leq 2+\frac{1}{x}-\ln 4.\Ee
\eprop
Comparing (\ref{eq:psiany}) and (\ref{eq:diffpsi1}) for $y=\frac 12$, we deduce the following.
\bcor
For $x>0$,
    \Be\label{eq:psismallfinal}
    \max\left\{1+\frac{1}{x}-\frac{1}{2x+1}-\ln 4,\frac{1}{x}-\ln\left(1+\frac{1}{2x}\right)\right\}\leq \Psi\left(x+\frac 12\right)-\Psi(x)\leq \frac{x+1}{x(2x+1)}.
    \Ee
\ecor
We note that the above result provides bounds for $0<x\leq 1$ that could be new in the literature. Unfortunately for $x>1$, it does not improve (\ref{eq:Qiu}). Note also that in the above discussion, we did not use (\ref{eq:diffquotient2}) because its bounds are not better than those of (\ref{eq:diffquotient200}).
\section{Conclusion}
In this paper, we have obtained several new inequalities for ratios of the beta function. We have shown our bounds improved the ones already available in the literature. We have also obtained some bounds for the beta function and the difference $\Psi(x+y)-\Psi(x)$.

\section{Compliance with Ethical Standards}
\begin{itemize}
\item {\bf Funding}

There is no funding support to declare.

\item {\bf Disclosure of potential conflicts of interest}

The authors have no relevant financial or non-financial interests to disclose.

\item {\bf Author Contributions}

All authors contributed to the conception and design of the study.  All authors read and approved the final manuscript.

\item {\bf Data availability statements}

Data sharing is not applicable to this article, as no data was created or analyzed in this study. 
\end{itemize}
\section{Appendix}
The following is the python code used to create the plot in Figure \ref{figlarge}. The code was written in Python 3.12.12 with Numpy version 2.3.4 and Matplotlib version 3.10.7.
%\appendix
\begin{lstlisting}[language=Python]
import numpy as np
import matplotlib.pyplot as plt
from matplotlib import colors
from mpl_toolkits.axes_grid1 import make_axes_locatable

def B(a,b,y):
    from numpy import log, exp
    c = 1 - (y/(b-a))*log((b+y)/(a+y))
    alpha = log(c) - 1
    beta = 1/c
    return exp((alpha+beta)*(b-a))*((a+y)/(b+y))**(y*beta)

def A(a,b,y):
    from numpy import log, exp
    c = 1 - (y/(b-a))*log((b+y+1)/(a+y+1))
    alpha = log(c) - 1
    beta = 1/c
    return exp((alpha+beta)*(b-a))*((a+y+1)/(b+y+1))**(y*beta)*(a/b)*((b+y)/(a+y))**y

F = lambda a,b,y: np.sign(A(a,b,y)-B(a,b,y))
a = np.linspace(0,500,2001,dtype=np.float128)
a[0] = 1e-16
n = a.size
tck_pos = list(np.linspace(0,a.size-1,11).astype(int))
tck_vals = [int(a[i]) for i in tck_pos]

y = np.array([1e-16, 0.2, 0.4, 0.6, 0.8, 1.0, 1.000001, 1.000005,
             1.000008, 1.00001, 1.00003, 1.00007, 1.0002, 1.2])

a,b = np.meshgrid(a,a)
a = a[:,:,np.newaxis]
b = b[:,:,np.newaxis]
b = np.flipud(b)
y = y[np.newaxis,np.newaxis,:]
Z = F(a,b,y)

fig, ax = plt.subplots(3,5,figsize=(12,8),sharex=True,sharey=True,constrained_layout=True)
cmap = plt.cm.grey
norm = colors.BoundaryNorm(np.array([-1.5,-0.5,0.5,1.5]), cmap.N)

for k in range(y.size):
    im = ax[k//5,k%5].imshow(Z[:,:,k],extent=(0,n,0,n),cmap=cmap,norm=norm)
    ax[k//5,k%5].set_xticks(tck_pos,tck_vals,rotation=90);
    ax[k//5,k%5].set_yticks(tck_pos,tck_vals);
    ax[k//5,k%5].set_title(r'$y='+str(y[0,0,k])+r'$')

fig.supxlabel('a');
fig.supylabel('b');
divider = make_axes_locatable(ax[-1,-1])
ax[-1,-1].axis('off')
cax = divider.append_axes('left', size='10%', pad=0.05)
cbar = fig.colorbar(im,ticks=np.array([-1,0,1]),cax=cax)
cbar.ax.set_yticklabels([r'$-1,\ B(a,b,y)>A(a,b,y)$',r'$0,\ B(a,b,y)=A(a,b,y)$',r'$1,\ A(a,b,y)>B(a,b,y)$']);
cbar.ax.set_title('Colour Key');
\end{lstlisting}

\bibliographystyle{plain}

\end{document}